\documentclass[a4paper,12pt]{amsart}
\usepackage[utf8]{inputenc}
\usepackage{hyperref}
\usepackage{fullpage}
\usepackage{amsmath, amssymb}
\usepackage[alphabetic,nobysame, initials]{amsrefs}
\usepackage{graphicx}

\newtheorem{thm}{Theorem}[section]
\newtheorem{lem}[thm]{Lemma}
\theoremstyle{definition}
\newtheorem{defn}[thm]{Definition}
\newtheorem{claim}[thm]{Claim}
\newtheorem{rem}[thm]{Remark}

\renewcommand{\Re}{\mathbb R}
\renewcommand{\epsilon}{\varepsilon}
\newcommand{\Red}{\Re^d}

\newcommand{\B}{\mathbf B}
\newcommand{\FF}{\mathcal F}
\newcommand{\GG}{\mathcal G}

\renewcommand{\phi}{\varphi}
\newcommand{\vol}[1]{\operatorname{vol}\left(#1\right)}

\newcommand{\card}[1]{\left|#1\right|}
\newcommand{\bd}{\operatorname{bd}}
\newcommand{\spann}{\operatorname{span}}
\newcommand{\tr}{\operatorname{tr}}
\newcommand{\conv}{\operatorname{conv}}
\newcommand{\iprod}[2]{\left< #1,#2\right>}

\title{Proof of a conjecture of B\'ar\'any, Katchalski and Pach}
\author[M. Nasz\'odi]{M\'arton Nasz\'odi}
\address{
M\'arton Nasz\'odi,
EPFL,
SB MATHGEOM DCG
CH-1015 Lausanne
Switzerland\\
\and
ELTE, 
Dept. of Geometry,
Lorand E\"otv\"os University,
P\'azm\'any P\'eter S\'et\'any 1/C
Budapest, Hungary 1117
}

\email{marton.naszodi@math.elte.hu}

\begin{document}
\begin{abstract}
In \cite{BKP84}, B\'ar\'any, Katchalski and Pach proved the following 
quantitative form of Helly's theorem. If the intersection of a family of convex 
sets in 
$\mathbb{R}^d$ is of volume one, then the intersection of some subfamily of at 
most 
$2d$ members is of volume at most some constant $v(d)$. In \cite{BKP82}, the 
bound $v(d)\leq 
d^{2d^2}$ is proved and $v(d)\leq d^{cd}$ is conjectured. We confirm it.
\end{abstract}

\maketitle

\section{Introduction and Preliminaries}

\begin{thm}\label{thm:qvol}
 Let $\FF$ be a family of convex sets in $\Red$ such that the volume of its 
intersection is $\vol{\cap\FF}>0$. 
Then there is a subfamily $\GG$ of $\FF$ with $\card\GG\leq 2d$ and
$\vol{\cap\GG}\leq Ce^{d}d^{2d} \vol{\cap\FF}$, where $C>0$ is a universal 
constant.
\end{thm}

We recall the note from \cite{BKP84} that the number $2d$ is optimal, as shown 
by the $2d$ half-spaces supporting the facets of the cube.

We introduce notations and tools that we will use in the proof.
We denote the closed unit ball centered at the origin $o$ in the 
$d$-dimensional Euclidean space $\Red$ by $\B$. For the scalar product of 
$u,v\in\Red$, we use $\iprod{u}{v}$, and the length of $u$ is 
$|u|=\sqrt{\iprod{u}{u}}$. The tensor product $u\otimes u$ is 
the rank one linear operator that maps any $x\in\Red$ to the vector 
$(u\otimes u)x=\iprod{u}{x}u\in\Red$. For a set $A\subset\Red$, we denote its 
polar by $A^{\ast}=\{y\in\Red:\iprod{x}{y}\leq 1,\mbox{ for all }x\in A\}$. The 
volume of a set is denoted by $\vol{\cdot}$.

\begin{defn}
 We say that a set of vectors $w_1,\ldots,w_m\in \Red$ with weights 
$c_1,\ldots,c_m>0$ form a \emph{John's decomposition of the identity}, if 
  \begin{equation}\label{eq:john}
    \sum_{i=1}^{m} c_iw_i = o
  \;\;\mbox{ and }\;\;
    \sum_{i=1}^{m} c_i w_i\otimes w_i = I,
  \end{equation}
where $I$ is the identity operator on $\Red$.
\end{defn}

A \emph{convex body} is a compact convex set in $\Red$ with non-empty interior.
We recall John's theorem \cite{Jo48} (see also \cite{Ba97}).
\begin{lem}[John's theorem]\label{lem:john}
  For any convex body $K$ in $\Red$, there is a unique ellipsoid of maximal 
  volume in $K$. Furthermore, this ellipsoid is $\B$ if, and only if, there are 
  points $w_1,\ldots,w_m\in\bd{\B}\cap\bd{K}$ (called \emph{contact points}) 
and corresponding weights 
$c_1,\ldots,c_m>0$ that form a John's decomposition of the identity.
\end{lem}

It is not difficult to see that if $w_1,\ldots,w_m\in\bd{\B}$ and corresponding 
weights $c_1,\ldots,c_m>0$ form a John's decomposition of the identity, then 
$\{w_1,\ldots,w_m\}^\ast\subset d\B$, cf. \cite{Ba97} or Theorem 5.1 
in \cite{GLMP04}. By polarity, we also obtain that 
$\frac{1}{d}\B\subset\conv(\{w_1,\ldots,w_m\})$.

One can verify that if $\Delta$ is a simplex in $\Red$, and $E$ is the largest 
volume ellipsoid in $\Delta$, then
\begin{equation}\label{eq:simplexvr}
 \frac{\vol{E}}{\vol{\Delta}}=\frac{d!\vol{\B}}{d^{d/2}(d+1)^{(d+1)/2}}.
\end{equation}

We will use the following form of the Dvoretzky-Rogers lemma \cite{DvoRo50}.
\begin{lem}[Dvoretzky-Rogers lemma]\label{lem:dvoro}
 Assume that $w_1,\ldots,w_m\in\bd{\B}$ and $c_1,\ldots,c_m>0$ form a 
John's decomposition of the identity. Then there is an orthonormal basis 
$z_1,\ldots,z_d$ 
of $\Red$, and a subset $\{v_1,\ldots,v_d\}$ of $\{w_1,\ldots,w_m\}$ such that 
\begin{equation}\label{eq:dvoro}
  v_i\in\spann\{z_1,\ldots,z_i\},
  \;\;\mbox{ and }\;\;
  \sqrt{\frac{d-i+1}{d}}\leq \iprod{v_i}{z_i}\leq 1,
  \;\;\mbox{ for all }i=1,\ldots,d.\;\;
\end{equation}
\end{lem}
This lemma is usually stated in the setting of John's theorem, that is, when 
the vectors are contact points of a convex body $K$ with its maximal volume 
ellipsoid, which is $\B$. And often, it is assumed in the statement that $K$ is 
symmetric about the origin, see for example \cite{GiEtAl14}. Since we make no 
such assumption (in fact, we make no reference to $K$ in the statement of 
Lemma~\ref{lem:dvoro}), we give a proof in Section~\ref{sec:dvoro}.

\section{Proof of Theorem~\ref{thm:qvol}}

\begin{figure}[tb]
    \centering
    \includegraphics[width=0.3\textwidth]{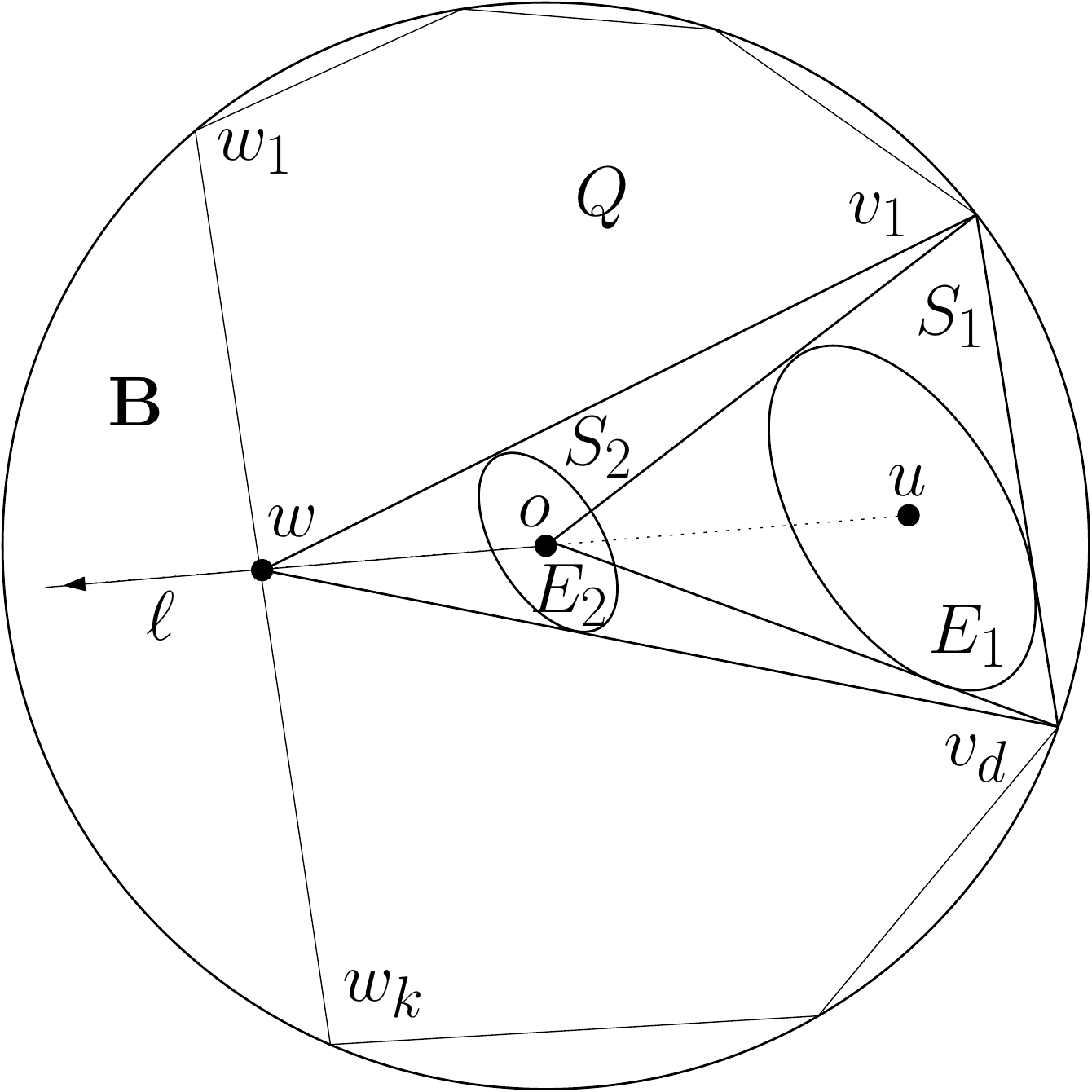}
    \caption[]{}
    \label{fig:bkp}
\end{figure}

Without loss of generality, we may assume that $\FF$ consists of closed 
half-spaces, and also that $\vol{\cap\FF}<\infty$, that is, $\cap\FF$ is a 
convex body in $\Red$. As shown in \cites{BKP84}, by continuity, we may 
also assume that $\FF$ is a finite family, that is $P=\cap\FF$ is a 
$d$-dimensional polyhedron.

The problem is clearly affine invariant, so we may 
assume that $\B\subset P$ is the ellipsoid of maximal volume in $P$. 

By Lemma~\ref{lem:john}, there are contact points 
$w_1,\ldots,w_m\in\bd{\B}\cap\bd{P}$ (and weights $c_1,\ldots,c_m>0$) that form 
a John's decomposition of the 
identity. We denote their convex hull by $Q=\conv\{w_1.\ldots,w_m\}$.
Lemma~\ref{lem:dvoro} yields that there is an orthonormal basis 
$z_1,\ldots,z_d$ 
of $\Red$, and a subset $\{v_1,\ldots,v_d\}$ of the contact points 
$\{w_1,\ldots,w_m\}$ such that \eqref{eq:dvoro} holds.

Let $S_1=\conv\{o,v_1,v_2,\ldots,v_d\}$ be the simplex spanned by these contact 
points, and let $E_1$ be the largest volume ellipsoid contained in $S_1$. We 
denote the 
center of $E_1$ by $u$. Let $\ell$ be the ray emanating from the origin in the 
direction of the vector $-u$. Clearly, the origin is in the interior of $Q$. In 
fact, by the remark following Lemma~\ref{lem:john}, $\frac{1}{d}\B\subset Q$. 
Let $w$ be the 
point of intersection of the ray $\ell$ with $\bd Q$. Then $|w|\geq 1/d$. Let 
$S_2$ denote the simplex $S_1=\conv\{w,v_1,v_2,\ldots,v_d\}$. See 
Figure~\ref{fig:bkp}.

We apply a contraction with center $w$ and ratio $\lambda=\frac{|w|}{|w-u|}$ on 
$E_1$ to obtain the ellipsoid $E_2$. Clearly, $E_2$ is centered at the origin 
and is contained in $S_2$. Furthermore,
\begin{equation}\label{eq:lambdanotsmall}
 \lambda=\frac{|w|}{|u|+|w|}\geq\frac{|w|}{1+|w|}\geq\frac{1}{d+1}.
\end{equation}

Since $w$ is on $\bd Q$, by Caratheodory's theorem, $w$ is in the convex hull 
of some set of at most $d$ vertices of $Q$. By re-indexing the vertices, we may 
assume that 
$w\in\conv\{w_1,\ldots,w_k\}$ with $k\leq d$. Now,
\begin{equation}\label{eq:ellipsoidcontained}
 E_2\subset S_2\subset \conv\{w_1,\ldots,w_k, v_1,\ldots,v_d\}.
\end{equation}

Let $X=\{w_1,\ldots,w_k, v_1,\ldots,v_d\}$ be the set of these
unit vectors, and let $\GG$ denote the family of those half-space 
which support $\B$ at the points of $X$. Clearly, $|\GG|\leq 2d$. Since the 
points of $X$ are contact points of $P$ and $\B$, we have that 
$\GG\subseteq\FF$. By \eqref{eq:ellipsoidcontained},
\begin{equation}\label{eq:GGcontained}
\cap\GG=X^{\ast}\subset E_2^{\ast}. 
\end{equation}

Since $\B\subset\cap\FF$, by \eqref{eq:GGcontained} and 
\eqref{eq:lambdanotsmall}, and \eqref{eq:simplexvr} we have
\begin{equation}\label{eq:finalestimate}
 \frac{\vol{\cap\GG}}{\vol{\cap\FF}}\leq
 \frac{\vol{E_2^{\ast}}}{\vol{\B}}
 =\frac{\vol{\B}}{\vol{E_2}}
 \leq(d+1)^d\frac{\vol{\B}}{\vol{E_1}}
 =
 \frac{d^{d/2}(d+1)^{(3d+1)/2}}{d!\vol{S_1}}.
\end{equation}
By \eqref{eq:dvoro},
\begin{equation}\label{eq:simplexvol}
 \vol{S_1} \geq \frac{1}{d!}\cdot 
\frac{\sqrt{d!}}{d^{d/2}}=\frac{1}{\sqrt{d!}d^{d/2}},
\end{equation}
which, combined with \eqref{eq:finalestimate}, yields the desired result, 
finishing the proof of Theorem~\ref{thm:qvol}.

\begin{rem}
In the proof, in place of the Dvoretzky-Rogers lemma, we could select
the $d$ vectors $v_1,\ldots,v_d$ from the contact points randomly: 
picking $w_i$ with probability $c_i/d$ for $i=1,\ldots,m$, and repeating this 
picking independently $d$ times. Pivovarov proved (cf. Lemma~3 in \cite{Pi10}) 
that the 
expected volume of the random simplex $S_1$ obtained this way is the same as 
the 
right 
hand side in \eqref{eq:simplexvol}.
\end{rem}

\section{Proof of Lemma~\ref{lem:dvoro}}\label{sec:dvoro}

We follow the proof in \cite{GiEtAl14}. 
\begin{claim}\label{claim:trT}
Assume that $w_1,\ldots,w_m\in\bd{\B}$ 
and $c_1,\ldots,c_m>0$ form a John's decomposition of the identity. 
Then for any linear map $T:\Red\to\Red$ there is an $\ell\in\{1,\ldots,m\}$ 
such that
\begin{equation}
 \iprod{w_\ell}{Tw_\ell}\geq\frac{\tr T}{d},
\end{equation}
where $\tr T$ denotes the trace of $T$.
\end{claim}

For matrices $A,B\in\Re^{d\times d}$ we use $\iprod{A}{B}=\tr\left(AB^T\right)$ 
to denote their Frobenius product.

To prove the claim, we observe that
\begin{equation*}
 \frac{\tr T}{d}=\frac{1}{d}\iprod{T}{I}=\frac{1}{d}\sum_{i=1}^m 
c_i\iprod{T}{w_i\otimes w_i}
 =\frac{1}{d}\sum_{i=1}^m c_i\iprod{Tw_i}{w_i}.
\end{equation*}
Since $\sum_{i=1}^m c_i=d$, the right hand side is a weighted average of the 
values $\iprod{Tw_i}{w_i}$. Clearly, some value is at least the average, 
yielding Claim~\ref{claim:trT}.

We define $z_i$ and $v_i$ inductively. First, 
let $z_1=v_1=w_1$. Assume that, for some $k<d$, we have found $z_i$ and $v_i$, 
for all $i=1,\ldots,k$. Let $F=\spann\{z_1,\ldots,z_k\}$, and let $T$ be 
the orthogonal projection onto the orthogonal complement $F^{\perp}$ of $F$. 
Clearly, $\tr T=\dim F^{\perp}=d-k$. By 
Claim~\ref{claim:trT}, for some $\ell\in\{1,\ldots,m\}$ we have
\begin{equation*}
 |Tw_\ell|^2=\iprod{Tw_\ell}{w_\ell}\geq\frac{d-k}{d}.
\end{equation*}
Let $v_{k+1}=w_\ell$ and $z_{k+1}=\frac{Tw_\ell}{|Tw_\ell|}$. Clearly, 
$v_{k+1}\in\spann\{z_1,\ldots,z_{k+1}\}$. Moreover,
\begin{equation*}
\iprod{v_{k+1}}{z_{k+1}}=\frac{\iprod{Tw_\ell}{w_\ell} }{|Tw_\ell|}=
\frac{|Tw_\ell|^2}{|Tw_\ell|}=|Tw_\ell|\geq\sqrt{\frac{d-k}{d}},
\end{equation*}
finishing the proof of Lemma~\ref{lem:dvoro}.

Note that in this proof, we did not use the fact that, in a John's 
decomposition 
of the identity, the vectors are balanced, that is $\sum_{i=1}^m c_iw_i=o$.

\bibliographystyle{amsalpha}

\end{document}